\documentclass[twocolumn,final,10pt]{amsart}
\pdfoutput=1
\usepackage{geometry}                
\usepackage[parfill]{parskip}    
\usepackage{graphicx}
\usepackage{amssymb, latexsym}
\usepackage{epstopdf}
\usepackage{gensymb}			

\usepackage{hyperref}

\begin{document}
\title{Paper On Statistics Of Repetitions}
\author{Alan M. Turing}

\begin{abstract}
This is a typeset version of Alan Turing's declassified Second World War paper \textit{Paper on Statistics of Repetitions.} See the companion paper \textit{The Applications of Probability to Cryptography}, available at \url{www.arxiv.org/abs/arXiv:1505.04714}, for Editor's Notes.
\vspace{0.2 in}
\newline
\textsc{Copyright.} The underlying manuscript is held by the National Archives in the UK  and can be accessed at \url{www.nationalarchives.gov.uk} using reference number HW 25/38. Readers are encouraged to obtain a copy. 
The original work was under Crown copyright, which has now expired; the work is now in the public domain.
\end{abstract}

\maketitle
%
%
\begin{flushleft}
	STATISTICS OF REPETITIONS
\end{flushleft}
In order to be able to obtain reliable estimates of the value of given repeats we need to have information about repetition in plain language. Suppose for example that we have placed two messages together and that we find repetitions consisting of a tetragramme, two bigrammes, and fifteen single letters, and that the total \textit{overlap} was 105, \textit{i.e.} that the maximum possible number of repetitions which could be obtained by altering letters of the messages is 105; suppose also that the lengths of the messages are 200 and 250; in such a case what is the probability of the fit being right, no other information about the day's traffic being taken into consideration, but information about the character of the unenciphered text being available in considerable quantity?

In theory this can be solved as follows. We take a vast number of typical decodes, say $10^{10^{6}}$, and from them we select all of length 200 and all of length 250. We encipher all of these messages at all possible positions on the machine (neglecting for simplicity the complications due to different daily keys). We then compare each message 200 long with each 250 long in such a way as to get an overlap of 105 as with the fit under consideration. From the resulting comparisons we pick out just those cases where the repetitions have precisely the same form as in the case in question. 

This set of comparisons will be called the \textit{relevant} comparisons. Among the relevant comparisons there will be some which are \textit{right} comparisons, \textit{i.e.} where corresponding letters of the two messages were enciphered with the same position of the machine. The probability that our original fit was right can now be expressed in the form:
\[
	\frac{\textit{Number of right relevant comparisons}}{\textit{Total number of relevant comparisons}}.
\]
%
%
The work involved in this theoretical method can be vastly reduced if we make a few harmless assumptions. In the first place if we assume that the encipherment keys at the various positions of the machine are \textit{hatted} we can calculate the number of relevant wrong comparisons. Suppose the total number of repeated letters in the case in question is $R$, then 
\begin{gather*}
	\frac{\textit{Number of relevant wrong comparisons}}{\textit{Total number of wrong comparisons}} \\
	= \left( \frac{L}{26} \right)^{R} \left( \frac{25}{26!} \right)^{\left( L - R \right)}.
\end{gather*}
For the calculation of the number of relevant right comparisons we have to make other assumptions. The sort of assumption that we need is that a repetition in one place is not made any the more or less likely by a repetition elsewhere. Actually this assumption would not be quite true, as it clearly does not hold in the case of adjacent letters. For more practical purposes I think the following assumption is sufficiently near to the truth:
 
 \textit{If we know that at a certain point P there is not a repetition, then knowledge that there is or is not a repetition at a point A before P does not make a repetition at a point B after P either more likely or less likely. The probability of a repetition at any point is also independent of its distance from the end of either message.}
 
With these assumptions we could get the right distribution of numbers of comparisons between the various repetition figures if we assume the repetition figures for the comparisons constructed in this way. We are given an urn containing a large number of cards, some bearing the words \textbf{no repeat}, some bearing \textbf{simple repeat}, some \linebreak\textbf{bigramme}, some \textbf{trigramme}, and so on. To construct a random sample of repetition figures of comparisons of given length we make a series of draws from the urn.
 
%
%
The first few draws determine the repetition figure for the first comparison, the next few for the next comparison, and so on. When we draw \textbf{no repeat} we have to add a \texttt{`O'} to the repetition figure, when we draw \textbf{simple repeat} add \texttt{`XO'}, for \textbf{bigramme} we add \texttt{`XXO'} and so on. When we have got to the right length of overlap required the comparison is completed and our next draws refer to the next comparison. If it happens that the right length is never reached because we `jump past it' then we scrap that comparison, and go on to the next. As an example suppose that we are making comparisons with an overlap of 12, and that our first draws are \textbf{tetragramme},  \textbf{no rep},  \textbf{no rep},  \textbf{no rep},   \textbf{bigramme},  \textbf{no rep},  \textbf{trigramme},  \textbf{13-gramme}, then  \textbf{no rep} 13 times, our first two comparisons will have the repetition figures:
 \texttt{
 \begin{center}
 	XXXXOOOOXXOO \\
	OOOOOOOOOOOO
 \end{center}}
 the one starting \texttt{`XXXO'} being rejected because we never reach the right length of overlap. (This arrangement requires that every repetition figure should end with \texttt{`O'}, and therefore the genuine repetition figure should be obtained by crossing this off; but I shall not be too meticulous about details arising from the ends of the comparison). 
 
The number of draws required to produce a given figure is the number of non repeating letters, \textit{i.e.} the overlap less the number of repeating letters. With our convention about crossing off the last letter we have to add 1.
 
 Two problems arise from this picture
 \begin{enumerate}
 	\item How do we calculate the correct proportions of cards in the urn?
 	\item Given the proportion of the cards in the urn, how do we calculate the number of right relevant comparisons, and hence the probability of a given fit?
 \end{enumerate}
 
%
%
The correct proportion of the cards in the urn can be calculated from the actual distribution of repetitions in the case of messages correctly set, or, what comes to the same thing, in messages unenciphered and arbitrarily set. Let us suppose that we have a large number of such comparisons of unenciphered messages, and that the messages are sufficiently long that complications arising from the ends of the messages can be neglected. 

	The proportion of cards bearing the words  \textbf{simple repeat},  \textbf{bigramme},  \textbf{trigramme}, etc., must obviously be in the same ratio as the number of corresponding repeats in our comparisons. The number of \textbf{no repeat} cards will be calculated slightly differently as we have to subtract one case of  \textbf{no repeat} for each sequence of repeating letters.
 
To get the best value from given material we naturally make every possible comparison. If we do this the right number of repetitions can be calculated quite easily without actually making the comparisons. Theoretically we can imagine the complete set of comparisons made in this way. 

First of all we write out all the decodes (say 50 of them) one after another round a circle; suppose that the number of letters on this circle is~$N$. The whole is then repeated on a concentric circle. All possible comparisons can be made by rotating the one circle with respect to the other. From these we have to remove the comparison in which the circles are not rotated at all, for obvious reasons. Also when the rotation is more than 180\degree \, we get essentially the same comparison as one with less than 180\degree. The net effect of this, taking into account also the special case of exact 180\degree \, rotation, is that the total overlap of all the comparisons is:
\[
	\frac{N(N - 1)}{2}.
\]
%
%
	Now let us consider for example the total number of tetragramme repeats in all these comparisons. These can be divided into repeats arising from \texttt{AAAA} those from \texttt{AAAB} $\dots$ those from \texttt{ZZZZ}, the largest contribution arising presumably from such tetragrammes as \texttt{EINS}. The number of tetragrammes arising from \texttt{EINS} consists of the number of \textit{pairs of hexagrams} such as \texttt{QEINSR}, \texttt{VEINSW} in which the first letters of each are different, the last different, and the remainder spell \texttt{EINS}. This number of pairs we will call the \textit{actual number of tetragramme repeats} arising from $EINS$. 

The \textit{actual number of tetragramme repeats} is obtained by summing over \texttt{AAAA}, \texttt{AAAB},  \dots  , \texttt{EINS},  \dots  , \texttt{ZZZZ}. This \textit{actual} number is not easily calculated directly, but we can more easily obtain the \textit{apparent number of tetragramme repeats}, and this leads to the actual number. The \textit{apparent number of tetragramme repeats arising from} \texttt{EINS} is defined to be the number of pairs of occurrences of \texttt{EINS} in the material, and the apparent number of tetragramme repeats defined by summation. 

We can also define the apparent number of tetragramme repeats in a comparison as the number of different series \texttt{XXXX} in the comparison. Thus a heptagramme repeats gives four apparent tetragramme repeats.The actual number of repeats can be calculated from the apparent in this way. 

Let $M_{r}$ be the apparent number of $r$-grammes, and  $N_{r}$ the actual number. 
 
Then
\[
	M_{r} =  N_{r} + 2N_{r+1} + 3N_{r+2} + \dotsc ,
\]
so that
\[
	M_{r} - M_{r+1} = N_{r} + N_{r+1} + N_{r+2} + \dotsc ,
\]
and
\begin{align*}
	N_{r} &= (M_{r} - M_{r+1}) - (M_{r+1} - M_{r+2}), \\
	&= M_{r} - 2M_{r+1} + M_{r+2} \ .
\end{align*}
%
%
It is therefore sufficient to calculate only apparent numbers and to carry these two stages further that we want to go with the actual numbers. In practice octagramme repeats are so certain to be right that it will be sufficient to have statistics only as far a heptagrammes. We therefore need statistics of apparent numbers of repeats as far a 9-grammes. To get these numbers of apparent repeats is sufficient to take all the 9-grammes in the material (\textit{i.e.} on the circle) and to put them into alphabetical order. This can be done very conveniently by Hollerith. 

The number of trigramme repeats say can then be found very simply (although with a good deal of labour) by considering only the first three letters of each 9-gramme. Suppose we denote by  $t$ a typical trigramme and by $n_{t}$ the number of its occurrences, then the apparent number of trigramme repeats is 
\[
	\sum_t \frac{n_{t} (n_{t} - 1)}{2} \ .
\]
When calculating the proportion of cards in the urn we must remember that the total number of cards is not 
\[
	\frac{N(N - 1)}{2} \ ,
\]
but is less than this by 
\[
	\sum r N_{r} \ .
\]
In our later calculations it is convenient to regard the comparisons in wrong places as also constructed by drawing from an urn. In this case we easily see that the apparent number of $r$-grammes is 
\[
	\left( \frac{N(N - 1)}{2} \right) \left( \frac{1}{26^{r}} \right),
\]
and from this we deduce that the actual proportion of $r$-gramme cards is
\[
	\left( \frac{25}{26^{(r+1)}} \right)  \text{ and of \textbf{no repeat} cards is } \left( \frac{25}{26} \right).
\]

%
%
We now turn to the problem of calculating the probability of a given fit when we know the proportion $\alpha_{r}$ of $r$-gramme cards that are in the urn for each $r$. The calculation is going to be slightly complicated by the convention which we introduced, that not all drawings can lead to a comparison. We have therefore to calculate the proportion of draws which do lead to a comparison, \textit{i.e.} in which the length does not overshoot the mark. The answer is that as the length of the overlap tends to infinity the proportions tends to
\[
	\frac{1}{1 + \sum r \alpha_{r}} \ ,
\]
in the case of hatted material this is $ \left( 25/26 \right)$.

Now we put $A = 1 - \sum \alpha_{r}$. Consider a repetition figure in which there are $k_{r}$ $r$-grammes. Let the overlap be $L$. The number of \textbf{no repeat} cards drawn is $L + 1 - \sum (r +1) k_{r}$. The proportion of right draws which are relevant is
\[
	A^{L + 1 - \sum (r + 1) k_{r}} \prod_{r=1}^\infty \alpha_{r}^{k_{r}} \ ,
\]
and then the proportion of the right \textit{comparisons} which are relevant is (assuming $L$ reasonably large)
\[
	\biggl(1 + \sum r \alpha_{r} \biggr) A^{L + 1 - \sum \left(r + 1\right) k_{r}} \prod_{r=1}^\infty \alpha_{r}^{k_{r}}.
\]
Similarly calculating with the urn whose proportions were made up from hatted materiel we find for the proportion of wrong comparisons which are relevant
\[
	\left( \frac{26}{25} \right) \times \left( \frac{25}{26} \right)^{L + 1 - \sum (r + 1) k_{r}}  
		\times \prod_{r=1}^\infty \left( \frac{25}{26^{\left( r + 1 \right)}} \right) ^{k_{r}}.
\]

\renewcommand{\thefootnote}
 {\ensuremath{\fnsymbol{footnote}}}
 
Hence the odds\footnote{ The odds on an event are defined to be the probability of the event divided by the probability of its negation.} on our fit are
\begin{multline*}
	q = \lambda \frac{ 25(1 + \sum r \alpha_{r}) }{26} \left( \frac{26A}{25} \right) ^{L + 1 - \sum (r + 1) k_{r}} \\
	 	 \times \prod_{r=1}^{\infty} \left( \frac{26^{\left( r + 1 \right)} \alpha_{r}}{25} \right) ^{k_{r}},
\end{multline*}
where $\lambda$ is the \textit{a priori} odds. This is most conveniently written as:
\begin{multline*}
	\log q =  \log \lambda + \sum \mu_r{} k_{r} - \nu L \\
		  + \log \left[ \left(1 - \sum \alpha_{r} \right) \left( 1 + \sum r \alpha_{r} \right) \right],
\end{multline*}
where
\[
	\mu_{r} =  \log \left(  \frac{\alpha_{r}26^{\left( r + 1 \right)}}{25} \right) - \left(r + 1\right) \log \left( \frac{26A}{25} \right),
\]
and
\[
	\nu = \log \left(  \frac{25}{26A} \right) \fallingdotseq \sum \alpha_{r} - \frac{2}{51}.
\]

%
%
In the case of overlap zero there is a discrepancy of 
\[
	\log \bigg[ \left(1 - \sum \alpha_{r} \right) \left(1 + \sum r  \alpha_{r} \right) \bigg],
\]
due to the overlap not being long. This term is in any case microscopic.

\vspace{18pt}
\hrule
\begin{center}
	$\infty$
\end{center}
\vspace{6pt}
\hrule

\end{document}